\documentclass[11pt]{amsart}

\usepackage[all]{xy}
\allowdisplaybreaks

\usepackage{ifthen}
\usepackage{soul}
\usepackage{latexsym}
\usepackage{amsthm,amssymb}
\usepackage{amsbsy,amsfonts,amsmath}

\usepackage{a4}

\usepackage[ansinew]{inputenc}
\usepackage{graphicx,fancyhdr}

\usepackage[dvips, dvipsnames, usenames]{color}

\newtheorem{theorem}{Theorem}[section]

\newtheorem{lemma}[theorem]{Lemma}
\newtheorem{coro}[theorem]{Corollary}

\newtheorem{prop}[theorem]{Proposition}

\theoremstyle{definition}
\newtheorem{definition}[theorem]{Definition}
\newtheorem{example}[theorem]{Example}

\theoremstyle{remark}
\newtheorem{remark}[theorem]{Remark}

\newcommand{\vi}{\textbf{(i)} }
\newcommand{\vii}{\textbf{(ii)} }
\newcommand{\viii}{\textbf{(iii)} }
\newcommand{\viv}{\textbf{(iv)} }

\newcommand\Z{\mathbb Z}
\newcommand\N{\mathbb N}
\newcommand\I{\mathbb I}

\def\cN{\mathcal{N}}

\def\cX{\mathcal{X}}

\def\cC{\mathcal{C}}
\def\cW{\mathcal{W}}
\def\cR{\mathcal{R}}

\def\cQ{\mathcal{Q}}

\newcommand\re{\operatorname{re}}

\newcommand\Hom{\operatorname{Hom}}
\newcommand\End{\operatorname{End}}
\newcommand\Aut{\operatorname{Aut}}
\newcommand\id{\operatorname{id}}

\newcommand\Sb{\mathbb S}

\newcommand{\bm}{\mathbf{m} }
\newcommand{\bM}{\mathbf{M} }
\newcommand{\vtxgpd}{\bullet}

\def\cWt{\widetilde{\mathcal{W}}}
\def\st{\widetilde{\sigma}}
\def\et{\widetilde{\er}}

\def\nb{\mathbf{n}}
\def\eb{\mathbf{e}}

\def\nr{\mathrm{n}}
\def\er{\mathrm{e}}

\newcommand{\al}{\alpha}

\newcommand{\bK}{{\mathbb{K}}}

\newcommand{\bZ}{{\mathbb{Z}}}

\def\pf{\begin{proof}}
\def\epf{\end{proof}}

\begin{document}

\title[Bruhat order and nil-Hecke algebras for Weyl groupoids]{Bruhat order and nil-Hecke algebras for Weyl groupoids}
\author[Angiono, Yamane]{Iv\'an Angiono}
\address{FaMAF-CIEM (CONICET), Universidad Nacional de C\'ordoba,
Medina A\-llen\-de s/n, Ciudad Universitaria (5000) C\' ordoba, Rep\'
ublica Argentina.}
\email{angiono@famaf.unc.edu.ar}
\author[]{Hiroyuki Yamane}
\address{University of Toyama,
Faculty of Science,
Gofuku 3190, Toyama-shi, Toyama 930-8555, Japan}
\email{hiroyuki@sci.u-toyama.ac.jp}

\thanks{\noindent 2010 \emph{Mathematics Subject Classification.}
16T05, Secondary 17B37. \newline The work of I. A. was partially supported by
CONICET, FONCyT-ANPCyT, Secyt (UNC). The work of H. Y. was partially supported by
Japan's Grand-in-Aid for Scientific Research (C) 25400040.}

\maketitle

\begin{abstract}
We introduce nil-Hecke algebras for Weyl groupoids. We describe a basis and
some properties of these algebras which lead to a notion of Bruhat order for Weyl groupoids.
\end{abstract}

\section{Introduction}

One of the remarkable points of contragredient Lie superalgebras, which makes a difference with Lie algebras,
is the existence of different matrices and parities on the generators whose corresponding Lie superalgebras are isomorphic. These algebras are related by odd reflections \cite{Ser-root systems} but then a problem arises:
how to manage these reflections together with the Weyl group (of the even part) at the same time.
By \cite{Ser-Ves} the Grothendieck ring of these Lie superalgebras is obtained as invariant of a \emph{super Weyl groupoid}.
The Weyl groupoid considered here seems to be related with their super Weyl groupoids.

A similar context appears for Nichols algebras of diagonal type, which are generalizations of (the positive part of)
quantized enveloping algebras and Frobenius-Lusztig kernels. The existence of Weyl groupoid for Nichols algebras was started in \cite{H-Weyl grp} for diagonal braidings and continued in \cite{AHS}
in a more general context. An axiomatic study was initiated in \cite{HY-weylgpd},
where a Matsumoto-Tits type theorem was given and a relation with Lie superalgebras was proposed.
Indeed every contragredient Lie superalgebra has a Weyl groupoid by \cite{AA-WGCLSandNA}, see also \cite{AYY} for characteristic 0 case.
Some braidings of diagonal type are related with contragredient Lie superalgebras over fields of
characteristic zero \cite{AAY}, but this relation can be extended to positive characteristic \cite{AA-WGCLSandNA} since a Weyl groupoid can be attached to any (finite-dimensional) contragredient Lie superalgebra.

This Weyl groupoid presents several applications. In \cite{HY-shapov} Heckenberger and the second author obtained
the factorization formula of the Shapovalov determinants of the generalized quantum groups
$U$, where an example of $U$ is any
Frobenius-Lusztig kernel (i.e., Lusztig's small quantum group); they used an action of
the Weyl groupoid on Verma modules of $U$.
Any reduced expression of its longest element gives an explicit realization of the universal R-matrix of $U$
\cite{AY2015}.
A special reduced expression of it was used to classify finite dimensional irreducible
representations of $U$ \cite{AYY}.
A weak order of the Weyl groupoids turned out to be useful since it is used in the classification
of the coideals of Nichols algebras \cite{HS}. These orders are related with convex orders on the positive roots \cite{A-convex}, a notion generalizing the one for classical root systems.
In \cite{HW} topological structures were studied, involving hyperplanes and associated with Weyl groupoids. It was used the left order of the Weyl groupoid.

Then one can naturally ask if there exists an analogous of Bruhat order for Weyl groupoids.

The Bruhat order of Coxeter groups essentially appears in various areas in mathematics,
for example in so-called Schubert calculus treated for studying the cohomology of a flag manifold \cite{BS}.
Another example is the formula given in \cite[Proposition 2.2]{IN}, a key result to concretely calculate
the value of the equivariant Schubert class.

In this paper, we introduce a Bruhat order, or a strong order,
for the elements of the Weyl groupoids. To this end we also introduce nil-Hecke algebras for Weyl groupoids and obtain some properties.
A future work is to establish the existence of Kazhdan-Lusztig polynomials of Weyl groupoids and the applications of these results to
the representation theory of Nichols algebras.

This paper is organized as follows. We give a brief introduction of Bruhat orders for Coxeter groups in Section \ref{section:bruhat Coxeter group}
including the main properties. Then in Section \ref{section:coxeter-groupoid} we recall the definition of basic datum, Coxeter groupoid and generalized root systems. In Section \ref{section:Nil hecke alg} we introduce the nil-Hecke algebras for Weyl groupoids.
We describe a basis indexed by the elements of the Weyl groupoid and a representation of a groupoid covering the Weyl groupoid on the nil-Hecke algebra.
Finally we introduce the definition of Bruhat order in Section \ref{section:Bruhat order} and prove the independence of the reduced expression.

\subsection*{Notation}
For $\theta\in\N$, let $\I_\theta=\{1,2,\dots,\theta\}$, or simply $\I$ if $\theta$ is clear from the context.
Let $\{\al_i\}_{i\in\I}$ be the  canonical basis of $\Z^I$.
Also, ${\mathbb{K}}$ will denote a commutative ring.

\section{Bruhat order for Coxeter groups}\label{section:bruhat Coxeter group}

We start by recalling the definition the Bruhat order of Coxeter groups
and we outline a proof about his good definition using Nil-Hecke algebras.

A \emph{Coxeter matrix} is a symmetric matrix $M = (m_{ij})_{i,j \in \I}$ with entries in
$\N \cup \{+\infty\}$ such that $m_{ii} = 1$ and $m_{ij} \geq 2$, for all $i\neq j\in \I$.
The \emph{Coxeter system} of $M$ is the pair $(W,S)$, where $W$ is the group presented by generators $S = \{s_1, \dots, s_{\theta}\}$ and relations
$(s_is_j)^{m_{ij}}= e$, for all $i, j\in \I$.

\begin{definition}
Let $(W,S)$ be a Coxeter system. Given $w$, $w^\prime\in W$, we say that
$w' < w$ if there exists a reduced expression $w=s_1s_2\cdots s_n$ and a subsequence $i_1,i_2,\ldots,i_r$
of $1,2,\ldots,n$ such that $w'=s_{i_1}s_{i_2}\cdots s_{i_r}$ is a reduced expression of $w'$. It defines a partial order on $W$ called the \emph{Bruhat order}.
\end{definition}

\begin{theorem}
Let $w$, $w'\in W$ be such that $w' < w$. Then for any reduced expression $w=t_1t_2\cdots t_n$, $t_k\in S$, there exists a subsequence
$i_1,i_2,\ldots,i_r$ of $1,2,\ldots,r$ such that $w'=t_{i_1}t_{i_2}\cdots t_{i_r}$ is a reduced expression of $w'$.
\end{theorem}

There exists a proof involving the nil-Hecke algebra of $(W,S)$ when $W$ is a Weyl group.
Let us explain it. Let ${\mathcal{N}}$ be a free left ${\mathbb{K}}$-module
with basis $n_w$, $w\in W$. Hence ${\mathcal{N}}$ is a ${\mathbb{K}}$-algebra with multiplication
$$
n_wn_{w^\prime}:=
\left\{
\begin{array}{ll}
n_{ww^\prime}& \quad\mbox{if $\ell(ww^\prime)=\ell(w)+\ell(w^\prime) $}  \\
0 & \quad\mbox{otherwise}
\end{array}\right.
$$
The algebra ${\mathcal{N}}$ is called the {\it{Nil-Hecke algebra}}.

Let $\Delta$ be a root system associated to $(W,S)$.
Let $\Pi=\{\alpha_i|i\in I\}$ be the basis of $\Delta$, $\Delta^+:=\Delta\cup \N_0\Pi$.
Then $\Delta=\Delta^+\cup-\Delta^+$.


Assume ${\mathbb{K}}$ is the free additive group with basis $\Pi$.
Let $n_i:=n_{s_i}\in{\mathcal{N}}$. For $\lambda\in{\mathbb{K}}$,
let $h_i(\lambda):=1+\lambda n_i$. Given $s_{i_1}\cdots s_{i_m}$ a reduced expression of $w\in W$, define $h(w)\in{\mathcal{N}}$ by
\begin{equation*}
h(w):=h_{i_1}(\alpha_{i_1})h_{i_2}(s_{i_1}(\alpha_{i_2}))\cdots h_{i_{\ell(w)}}(s_{i_1}\cdots s_{i_{\ell(w)-1}}(\alpha_{i_{\ell(w)}})).
\end{equation*}
It follows that $h(w)$ is independent of the choice of a reduced expression of $w$. This is proved by direct computation for rank-two Coxeter groups, and the general proof holds
since the defining relations involve only two letters of $S$. Hence the theorem follows since
$$\alpha_{i_1},s_{i_1}(\alpha_{i_2}),\ldots,s_{i_1}\cdots s_{i_{\ell(w)-1}}(\alpha_{i_{\ell(w)}})\in \N_0\Pi. $$

\section{Weyl groupoids and generalized root systems}\label{section:coxeter-groupoid}

We follow the notation of \cite{AA-WGCLSandNA}, see also \cite{CH,HY-weylgpd}.

\subsection{Basic data and Coxeter groupoids}
Given $\theta\in\N$, $\I=\I_\theta$, $\cX\neq\emptyset$ a non-empty set and $\rho: \I  \to \Sb_{\cX}$, the pair $(\cX, \rho)$ is a \emph{basic datum}
of base $\vert\cX\vert$ and size $\I$ if $\rho_i^2 = \id$ for all $i\in \I$.

We denote by $\cQ_{\rho}$ the quiver with vertices $\cX$ and arrows
\begin{align*}
  & \sigma_i^x := (x, i, \rho_i(x)), i\in \I, x\in \cX, & \mbox{with target }&t(\sigma_i^x) = x, \mbox{ source }s(\sigma_i^x) = \rho_i(x).
\end{align*}

The \emph{object change diagram} \cite{CH}, or simply the diagram, of $(\cX, \rho)$ is a simplified graphical description of
$\cQ_{\rho}$. It is the graph with bullets labeled with $\cX$ and one arrow between $x \neq y$, decorated with the label $i$,
for each pair $(x, i, y)$, $(y, i, x)$ such that $\rho_i(x) = y$. Therefore the loops at a vertex $x$ are deduced from the diagram
and $\theta$. $(\cX, \rho)$ is \emph{connected} if $\cQ_{\rho}$ is connected.

For any quotient of the free groupoid $F(\cQ_{\rho})$, $\sigma_{i_1}^x\sigma_{i_2}\cdots \sigma_{i_t}$ means
$$ \sigma_{i_1}^x\sigma_{i_2}^{\rho_{i_1}(x)}\cdots \sigma_{i_t}^{\rho_{i_{t-1}} \cdots \rho_{i_1}(x)},$$
omitting the implicit superscripts, uniquely determined to have compositions.

\medskip

Given a basic datum $(\cX, \rho)$,
a \emph{Coxeter datum} is a triple $(\cX, \rho, \bM)$, where $\bM = (\bm^x)_{x\in \cX}$,
$\bm^x = (m^x_{ij})_{i,j\in \I}$, is a family of Coxeter matrices such that
\begin{align}\label{eq:coxeter-datum}
s((\sigma^x_i\sigma_j)^{m^x_{ij}}) &= x,& \mbox{for all }i,j\in \I &\mbox{ and } x\in \cX.
\end{align}

The \emph{Coxeter groupoid} $\cW(\cX, \rho, \bM)$ \cite[Definition 1]{HY-weylgpd} is
the groupoid generated by the quiver $\cQ_{\rho}$ with relations
\begin{align}\label{eq:def-coxeter-gpd}
 (\sigma_i^x\sigma_j)^{m^x_{ij}} &= \id_x, &i, j\in \I,\,  &x\in \cX.
\end{align}
Notice that either $\sigma_i^x$ is an involution if $ \rho_i(x) = x$, or else $\sigma_i^x$ is the inverse arrow of
$\sigma_i^{\rho_i(x)}$ if $\rho_i(x) \neq x$.

\medskip

As in \cite{HY-weylgpd} we denote by $\cWt$ the corresponding quotient of the quiver groupoid $\cQ_\rho$
subject to the same relations as for $\cW$ except \eqref{eq:def-coxeter-gpd} for $i=j$; that is, we omit the relations
$\sigma^x_i\sigma^{r_i(x)}_i=\er^x$ for $i\in\I$, $x\in\cX$.
Then there exists a canonical projection $\bK\cWt\twoheadrightarrow \bK\cW$
such that $\st^x_i\mapsto \sigma_i^x$, $\et^x\mapsto \er^x$.

\begin{remark}\label{rmk:covering gpd}
Let $\bK$ be a commutative ring with 1. The groupoid algebra of $\cW(\cX, \rho, \bM)$ is generated by $\er^x$, $\sigma^x_i$, $x\in \cX$, $i\in \I$, with relations \eqref{eq:def-coxeter-gpd} and
\begin{align*}
\er^x\er^y&=\delta_{x,y} \er^x, & \sum_{x\in \cX}\er^x&=1, &
\sigma^x_i\er^{\rho_i(x)}&=\er^x\sigma^x_i=\sigma^x_i.
\end{align*}
\end{remark}

A \emph{covering} of a basic datum $(\cX, \rho)$ is a 3-uple $(\widetilde{\cX},\widetilde{\rho},F)$ such that
$(\widetilde{\cX},\widetilde{\rho})$ is a basic datum and $F:\widetilde{\cX}\to\cX$
is a surjective map such that $\rho_i\circ F=F\circ\widetilde{\rho_i}$ for all $i\in\I$.

Given two Coxeter data $(\cX, \rho, \bM)$, $(\widetilde{\cX}, \widetilde{\rho}, \widetilde{\bM})$ and $F:\cX'\to\cX$ such that
$(\widetilde{\cX},\widetilde{\rho},F)$ is a covering of $(\cX, \rho)$, we say that $F$ is a \emph{covering of Coxeter data} if
$M^{F(t_1)}=M^x$ for all $x\in\widetilde{\cX}$. In such a case $F$ induces a canonical surjective map of groupoids
$F:\cW(\widetilde{\cX}, \widetilde{\rho}, \widetilde{\bM})\twoheadrightarrow \cW(\cX, \rho, \bM)$. It is a particular case of a \emph{morphism} between Coxeter data \cite{AA-WGCLSandNA}.

\subsection{Generalized root systems and Weyl groupoids}\label{subsection:weyl-groupoid}

Recall that $C=(c_{ij})\in\Z^{\theta\times\theta}$ is a generalized Cartan matrix \cite{K} if for all $i\neq j\in\I$, $c_{ii}=2$, $c_{ij}\leq 0$,
and $c_{ij}=0$ if and only if $c_{ji}=0$.

A \emph{semi-Cartan graph} is a triple $(\cX, \rho, \cC)$  such that $(\cX, \rho)$ is a basic datum of size $\I$
and $\cC = (C^x)_{x\in \cX}$, $C^x = (c^x_{ij})_{i,j \in \I}$, $x\in \cX$,
is a bundle of generalized Cartan matrices satisfying
\begin{align}\label{eq:condicion Cartan scheme}
c^x_{ij}&=c^{\rho_i(x)}_{ij} &  \mbox{for all }&x \in \cX, \, i,j \in \I.
\end{align}
We define $s_i^x\in GL_{\theta}(\Z)$ by
\begin{align}\label{eq:reflection-x}
s_i^x(\alpha_j)&=\alpha_j-c_{ij}^x\alpha_i, & j&\in \I,&  i &\in \I, x \in \cX.
\end{align}
Then $s_i^x$ is the inverse of $s_i^{\rho_i(x)}$ by \eqref{eq:condicion Cartan scheme}.

\begin{definition}{\cite[Definition 1]{HY-weylgpd}}
Let $(\cX, \rho, \cC)$ be a semi-Cartan graph. A \emph{generalized root system} (for short a GRS) is a collection
$\cR:= \cR(\cX, \rho, \cC, \Delta)$,
where  $\Delta = (\Delta^x)_{x\in \cX}$ is a family of subsets $\Delta^x \subset \Z^{\I}$ such that for all $x \in \cX$ and all $i\neq j\in\I$,
\begin{itemize}
  \item $\Delta^x= \Delta^x_+ \cup \Delta^x_-$, where $\Delta^x_{\pm} := \pm(\Delta^x \cap \N_0^I) \subset \pm\N_0^I$,
  \item $\Delta^x \cap \Z \alpha_i = \{\pm \alpha_i \}$,
  \item $s_i^x(\Delta^x)=\Delta^{\rho_i(x)}$ ,
  \item $(\rho_i\rho_j)^{m_{ij}^x}(x)=x$, where $m_{ij}^x :=|\Delta^x \cap (\N_0\alpha_i+\N_0 \alpha_j)|$.
\end{itemize}
$\Delta^x_+$, respectively $ \Delta^x_-$, are called the set of \emph{positive}, respectively \emph{negative}, roots.
The \emph{Weyl groupoid} of $\cR$ is $\cW = \cW(\cX, \rho, \cC)$. It is called \emph{finite} if $R^x$ is finite for some
$x\in\cX$, or equivalently, for all $x\in\cX$; see \cite{CH} for other equivalences.
\end{definition}

\begin{remark}
  Let $(\cX, \rho, \cC)$ be a semi-Cartan graph with generalized root system $\cR$. Then the Weyl groupoid is a Coxeter groupoid \cite{HY-weylgpd}; more precisely, $\cW(\cX, \rho, \cC)=\cW(\cX, \rho, \bM)$, where $m_{ij}^x$ are defined as above for all $i,j\in\I$, $x\in\cX$.
\end{remark}

For any $x \in \cX$, the sets of \emph{real} roots at $x$ is
\begin{align}
\label{defrealroot} (\Delta^{\re})^x &= \bigcup_{y\in \cX}\{ w(\alpha_i): \ i \in \I, \ w \in \cW(y,x) \}.
\end{align}
Notice that $w(\Delta^x)= \Delta^y$ for all $w \in \cW(x, y)$. A semi-Cartan graph $(\cX, \rho, \cC)$ is \emph{Cartan}  if $\Delta^x= \bigcup_{y\in \cX}\{ w(\alpha_i): \ i \in \I, \ w \in \cW(y,x) \}$
is a finite GRS.

The notion of covering of Coxeter data $(\cX, \rho, \bM)$, $(\widetilde{\cX}, \widetilde{\rho}, \widetilde{\bM})$ extend canonically to
semi-Cartan graphs and GRS by
requiring, respectively, that $C^{F(x)}=C^x$ and $\Delta^{F(x)}=\Delta^x$ for all $x\in\widetilde{\cX}$.

\begin{example}\label{ex:finite-Nichols}
{\rm{
A semi-Cartan (or Cartan) graph $(\cX, \rho, \cC)$ is called \emph{standard} if all the matrices $C^x$, $x\in\cX$, are equal.
Let $(\cX, \rho, \cC)$ be a standard Cartan graph with root system $\cR$. Then $\cR$ is finite if and only if each
$C^x$ is of finite type, in which case each $\Delta^x$ is the root system of $C^x$.

Finite standard root systems with one point are Weyl groups, of type $A_2$, $B_2$ or $G_2$. Some other finite
standard Weyl groupoids are close to Lie superalgebras, and then with Nichols algebras of super type, see \cite{AAY}.
For example, the following Cartan graphs for $\theta=2$:
\begin{align}\label{eq:cartan-graph standard A2}
& \setlength{\unitlength}{1mm}
\begin{picture}(40,5)(0,0)
\put(0,0){\circle*{1.5}}\put(1,2){{\small{$A_2$}}}
\put(4, 0){\line(1,0){11}}\put(9,1){{\small{$2$}}}
\put(17,0){\circle*{1.5}}\put(18,2){{\small{$A_2$}}}
\put(21, 0){\line(1,0){11}}\put(26,1){{\small{$1$}}}
\put(34,0){\circle*{1.5}}\put(35,2){{\small{$A_2$}}}
\end{picture}  \\
& \setlength{\unitlength}{1mm}
\begin{picture}(40,5)(0,0)
\put(0,0){\circle*{1.5}}\put(1,2){{\small{$B_2$}}}
\put(4, 0){\line(1,0){11}}\put(9,1){{\small{$2$}}}
\put(17,0){\circle*{1.5}}\put(18,2){{\small{$B_2$}}}
\end{picture}
\label{eq:cartan-graph standard B2}
\end{align}
}}
\end{example}

\section{Nil Hecke algebras for Weyl groupoids of finite rank}\label{section:Nil hecke alg}

\subsection{Basic properties of nil Hecke algebras}

We assume that \emph{$\cX$ is finite}. We define a nil-Hecke algebra for each Weyl groupoid, which generalizes the corresponding notion for Weyl groups. Next we prove the first properties of this algebra. In particular we describe a basis using methods close to those for nil-Hecke algebras over Weyl groups.

\begin{definition}
Let $\bK$ be a commutative ring with 1. Let $(\cX, \rho, \cC)$ be a semi-Cartan graph and $\cR$ a GRS of $(\cX, \rho, \cC)$.
The \emph{nil-Hecke algebra} of $(\cX, \rho, \cC)$ is the $\bK$-algebra
$\cN=\cN(\cX, \rho, \cC)$ generated by
$\er^x$, $n_i^x$, $x\in\cX$, $i\in I$, with relations
\begin{align}
\er^x \er^y&=\delta_{x,y} \er^x, & \sum_{x\in\cX} \er^x&=1, \notag \\
n_i^x \er^{\rho_i(x)}&= \er^xn_i^x=n_i^x, & n_i^{\rho_i(x)}n_i^x&=0, \label{eq:definign relations N} \\
\underbrace{n_i^xn_jn_i \dots}_{m_{ij}^x} &=\underbrace{n_j^xn_in_j \dots}_{m_{ij}^x}, \notag
\end{align}
for all $x,y\in\cX$, $ i\neq j\in I$.
\end{definition}

\begin{remark}\label{rem:proj from Wtilde}
There exists a surjective map $\Upsilon:\bK\cWt\twoheadrightarrow\cN$ given by $\et^x\mapsto \er^x$, $\st_i^x\mapsto n_i^x$, for all $x\in\cX$, $i\in I$.
\end{remark}

We use a strategy close to the one in \cite[Section 4.1]{M} to produce a basis of this algebra.
Let $E^x, L_i^x:\bK\cW\to \bK\cW$, $x\in\cX$, $i\in I$, be the linear maps
\begin{align*}
E^x(w)&:=\left\{ \begin{array}{ll} w, & w\in\Hom(\cW,x) \\ 0, & \mbox{otherwise,}\end{array} \right. \\
L_i^x(w)&:=\left\{ \begin{array}{ll} \sigma_i^x w, & \sigma_i^x w\neq 0, \, \ell(\sigma_i^x w)>\ell(w) \\ 0, & \mbox{otherwise.}\end{array} \right.
\end{align*}

\begin{lemma}\label{lemma:representation nil Hecke over kW}
There exists an algebra map $\Lambda:\cN\to\End(\bK\cW)$ such that
\begin{align*}
\er^x&\mapsto E^x, &  n_i^x&\mapsto L_i^x, & x\in\cX, \, &i\in I.
\end{align*}
\end{lemma}
\pf
It follows that
\begin{align*}
\sum_{x\in\cX} E^x&=\id_V, & E^x E^y&=\delta_{x,y}E^x,  & L_i^x E^{\rho_i(x)}&=E^xL_i^x=L_i^x,
\end{align*}
for all $x,y\in\cX$, $i\in\I$. Also, if $w\in\Hom(\cW,x)$ is such that $\ell(\sigma_i^{\rho_i(x)} w)>\ell(w)$, then $\sigma_i^{\rho_i(x)} w\in\Hom(\cW,\rho_i(x))$ satisfies $\sigma_i^x \sigma_i^{\rho_i(x)} w =w$, so $L_i^xL_i^{\rho_i(x)}(w)=0$; if $\ell(\sigma_i^{\rho_i(x)} w)>\ell(w)$, then $L_i^{\rho_i(x)}(w)=0$ by definition. Hence $L_i^x L_i^{\rho_i(x)}=0$. As
$$ \underbrace{\sigma_i^x\sigma_j\sigma_i \dots}_{m_{ij}^x}=\underbrace{\sigma_j^x\sigma_i\sigma_j \dots}_{m_{ij}^x}, $$
we deduce that
$$ \underbrace{L_i^xL_jL_i \dots}_{m_{ij}^x}=\underbrace{L_j^xL_iL_j \dots}_{m_{ij}^x}. $$
Therefore $\Lambda$ is well-defined.
\epf

\begin{theorem}\label{thm:basis nil Hecke}
Let $w\in\Hom(y,x)$, $x,y\in\cX$ be an element of length $m$. If $w=\sigma_{i_1}^x \cdots \sigma_{i_m}= \sigma_{j_1}^x \cdots \sigma_{j_m}$ are two reduced expressions, then
\begin{equation}\label{eq:indep reduced expression nil hecke}
n_{i_1}^x n_{i_2} \cdots n_{i_m}=n_{j_1}^x n_{j_2} \cdots n_{j_m}.
\end{equation}
Call then $T_w$ to this element, and set $T_{\id_x}:=\er^x$ for all $x\in\cX$. Then
\begin{equation}\label{eq:basis nil hecke algebra}
\{ T_w | w\in\Hom(y,x),\, x,y\in\cX\}
\end{equation}
is a basis of $\cN$.
\end{theorem}
\pf
We claim that $\Phi:\cN\to \bK\cW$, $\Phi(n)=\Lambda(n)(1)$ is an isomorphism. Given an element $w\in\cW$ of length $m$, take a reduced expression $w=\sigma_{i_1}^x \cdots \sigma_{i_m}$. Then
\begin{align*}
\Phi\left(n_{i_1}^x n_{i_2} \cdots n_{i_m}\right)&=L_{i_1}^x L_{i_2} \cdots L_{i_m}(1)=\sigma_{i_1}^x \sigma_{i_2} \cdots \sigma_{i_m}=w.
\end{align*}
Therefore $\Phi$ is surjective.
\medskip

To prove the injectivity we fix a reduced expression $w=\sigma_{i_1}^x \cdots \sigma_{i_m}$ for each $w\in\cW$. We set
$T_w':=n_{i_1}^x n_{i_2} \cdots n_{i_m}$ if $\ell(w)\geq 1$, and $T_{\id_x}:=\er^x$ for each $x\in\cX$. Let $\cN'$ be the subspace spanned by $T_w'$, $w\in\cW$. Note that $1\in\cN'$. For each $x,y,z\in\cX$ and $w\in\Hom(y,z)$, it holds that $\er^x T_w'=\delta_{x,z} T_w'$. Given $i\in I$, we compute the product $n_i^x T_w'$. If $\ell(w)=0$, that is $w=\id_y$ for some $y\in\cX$, then $n_i^x T_w'=\delta_{x,y}n_i^x\in\cN'$. Otherwise, let $w=\sigma_{i_1}^y \cdots \sigma_{i_m}$ be the fixed reduced expression; we consider three possibilities:
\begin{itemize}
\item $y\neq \rho_i(x)$: $n_i^x T_w'=0$.
\item $y=\rho_i(x)$, $\ell(\sigma_i^x w)=m+1$: let $v=\sigma_i^x w$, and $v= \sigma_{j_0}^x \sigma_{j_1} \cdots \sigma_{j_m}$ the fixed reduced expression.
By \cite[Theorem 5]{HY-weylgpd}, we have that
$$ \st_{j_0}^x \st_{j_1} \cdots \st_{j_m}= \st_{i}^x \st_{i_1} \cdots \st_{i_m} $$
in $\cWt$, so applying the morphism $\Upsilon$ in Remark \ref{rem:proj from Wtilde},
$$ n_i^x T_w'= n_i^x n_{i_1} n_{i_2} \cdots n_{i_m}=T_v'\in\cN'. $$
\item $y=\rho_i(x)$, $\ell(\sigma_i^x w)<m+1$: By \cite[Corollary 6]{HY-weylgpd}, there exists $j_1,\dots,j_{m+1}\in \I$, $t\in\{1,\dots,m\}$ such that $j_t=j_{t+1}$ and
$$ \st_{j_1}^x \st_{j_2} \cdots \st_{j_{m+1}}= \st_{i}^x \st_{i_1} \cdots \st_{i_m} $$
in $\cWt$, so by applying the morphism $\Upsilon$ in Remark \ref{rem:proj from Wtilde},
$$ n_i^x T_w' = n_{j_1}^x n_{j_1} \cdots n_{j_{m+1}}=0. $$
\end{itemize}
Then $n_i^x T_w'\in\cN'$, so $\cN'$ is a left ideal, and we conclude that $\cN'=\cN$; that is, the elements $T_w'$, $w\in\cW$, span $\cN$. As $\Phi$ apply them to a basis of $\bK\cW$, $\Phi$ is an isomorphism.
\medskip

If we set a different reduced expression $w=\sigma_{j_1}^x \cdots \sigma_{j_m}$, then
$$\Phi\left(n_{j_1}^x n_{j_2} \cdots n_{j_m}\right)=w$$
by an analogous computation, so \eqref{eq:indep reduced expression nil hecke} follows.
\epf

\begin{coro}
The representation $\Lambda:\cN\to\End(\bK\cW)$ of the nil-Hecke algebra $\cN$ of Lemma \ref{lemma:representation nil Hecke over kW} is faithful.
\end{coro}
\pf
It follows since the isomorphism $\Phi$ introduced in the previous proof factorizes through $\Lambda$.
\epf

The existence of coverings induce algebra maps between the corresponding nil-Hecke algebras.

\begin{coro}
Let $(\cW',f)$ be a covering of a Weyl groupoid $\cW$. Let $\cN'$ be the nil Hecke algebra of $\cW'$. There exists an injective algebra map $F:\cN\to\cN'$ such that
\begin{align}
\er^x&\mapsto \sum_{y\in f^{-1}(t_1)} \er^y, & \nr^x_i&\mapsto \sum_{y\in f^{-1}(t_1)} \nr^y_i, & x\in\cX,&\,i\in I.
\end{align}
\end{coro}
\pf
Set $\eb^x=\sum_{y\in f^{-1}(t_1)} \er^y$, $\nb^x_i=\sum_{y\in f^{-1}(t_1)} \nr^y_i$. Then these elements of $\cN'$ satisfy the corresponding defining relations \eqref{eq:definign relations N} of $\cN'$ for all $x\in\cX$, $i\neq j\in I$, so $F$ is well-defined.

The injectivity is a consequence of Theorem \ref{thm:basis nil Hecke}, because the image of that basis of $\cN$ is a linearly independent set in $\cN'$ by the same Theorem.
\epf

\subsection{A useful representation of $\widetilde{\cW}$}

Set $\bK=\Z[t_i|i\in \I]$.
We identify the additive subgroup $\oplus_{i\in \I}\bZ t_i$ of $\bK$ with $\Z^\I=\oplus_{i \in \I}\bZ\al_i$
via the group isomorphism $\Z^I\to\oplus_{i\in I}\bZ t_i$
defined by $\al_i\mapsto t_i$, and using this isomorphism we can define
$s:\oplus_{i\in I}\bZ t_i \to \oplus_{i\in I}\bZ t_i$ for each $s\in\Aut(\Z^I)$. It can be extended to a unique algebra map $f_s:\bK\to\bK$ such that $f_s(t_i)=s(t_i)$ for all $i\in I$. In this way we obtain an action of $\Aut(\Z^I)$ on $\bK$ by algebra maps.

\medskip

For $x\in \cX$, $i\in I$ and $t\in\bK$, let
\begin{equation*}
h^x_i(t):=
\left\{\begin{array}{ll}
\er^x+t \nr^x_i & \quad\mbox{if $\rho_i(x)=x$}, \\
t \nr^x_i  &  \quad\mbox{if $\rho_i(x)\ne x$}.
\end{array}\right.
\end{equation*}

\begin{lemma}\label{lemma:formula thm-one point}
Let $x\in \cX$ be such that $\rho_1(x)=\rho_2(x)=x$ and $d=m_{12}^x$ is finite. For $n\in\N_0$ set
$i_{2n}=2$, $i_{2n+1}=1$, $j_{2n}=1$, $j_{2n+1}=2$. Then
\begin{align*}
h_{i_1}^x(t_{i_1}) & h_{i_2}^x(s_{i_1}^x\cdot t_{i_2})\cdots h_{i_d}^x(s_{i_1}^x\cdots s_{i_{d-1}}\cdot t_{i_d})= \\
&=h_{j_1}^x(t_{j_1})h_{j_2}^x(s_{j_1}^x\cdot t_{j_2})\cdots h_{j_d}^x(s_{j_1}^x\cdots s_{j_{d-1}}\cdot t_{j_d}).
\end{align*}
\end{lemma}
\pf
By hypothesis $\Delta^x\cap\{\Z\alpha_1+\Z\alpha_2\}$ is a finite root system, so it corresponds to a finite Cartan matrix $C=\left( \begin{array}{cc}2 & c_{12}^x \\ c_{21}^x & 2 \end{array} \right)$. We have four possibilities.

\noindent
\vi \emph{$C$ is of type $A_1 \times A_1$}. Notice that
\begin{align*}
h^x_1(t_1)h^x_2(t_2)&= \er^x+ t_1 \, n_1^x +t_2 \, n_2^x+t_1t_2 \, n_1^xn_2^x \\
&= \er^x+ t_1 \, n_1^x+t_2 \, n_2^x + t_1t_2 \, n_2^xn_1^x = h^x_2(t_2)h^x_1(t_1).
\end{align*}

\noindent
\vii \emph{$C$ is of type $A_2$}. Then we compute
\begin{align*}
h^x_1(t_1)&h^x_2(t_1+t_2)h^x_1(t_2) = (\er^x+t_1 \nr^x_1)(\er^x+(t_1+t_2)\nr^x_2)(\er^x+t_2\nr^x_1) \\
& = \er^x+(t_1+t_2)\nr^x_1+(t_1+t_2)\nr^x_2+t_1(t_1+t_2)\nr^x_1\nr^x_2 \\
& \qquad +(t_1+t_2)t_2\nr^x_2\nr^x_1+t_1(t_1+t_2)t_2\nr^x_1\nr^x_2\nr^x_1 \\
& = (\er^x+t_2 \nr^x_2)(\er^x+(t_1+t_2)\nr^x_1)(\er^x+t_1 \nr^x_2)= h^x_2(t_2)h^x_1(t_1+t_2)h^x_2(t_1).
\end{align*}

\noindent
\viii \emph{$C$ is of type $B_2$}. Then we compute
\begin{align*}
h^x_1(t_1)&h^x_2(2t_1+t_2)h^x_1(t_1+t_2)h^x_2(t_2) = \\
& = (\er^x+t_1\nr^x_1)(\er^x+(2t_1+t_2)\nr^x_2)(\er^x+(t_1+t_2)\nr^x_1)(\er^x+t_2\nr^x_2)\\
& = \er^x+(2t_1+t_2)\nr^x_1+2(t_1+t_2)\nr^x_2+(2t_1+t_2)(t_1+t_2)\nr^x_1\nr^x_2 \\
& \qquad + (2t_1+t_2)(t_1+t_2)\nr^x_2\nr^x_1+t_1(2t_1+t_2)(t_1+t_2)\nr^x_1\nr^x_2\nr^x_1 \\
& \qquad + (2t_1+t_2)(t_1+t_2)t_2\nr^x_2\nr^x_1\nr^x_2+t_1(2t_1+t_2)(t_1+t_2)t_2 \nr^x_1\nr^x_2\nr^x_1\nr^x_2 \\
& = (\er^x+t_2\nr^x_2)(\er^x+(t_1+t_2)\nr^x_1)(\er^x+(2t_1+t_2)\nr^x_2)(\er^x+t_1\nr^x_1) \\
& = h^x_2(t_2)h^x_1(t_1+t_2)h^x_2(2t_1+t_2)h^x_1(t_1).
\end{align*}

\noindent
\viv \emph{$C$ is of type $G_2$}. Then we compute
\begin{align*}
h^x_1(t_1)&h^x_2(3t_1+t_2)h^x_1(2t_1+t_2)h^x_2(3t_1+2t_2)h^x_1(t_1+t_2)h^x_2(t_2) \\
& = 1+2(2t_1+t_2)\nr^x_1+2(3t_1+2t_2)\nr^x_2 \\
& \qquad +2(2t_1+t_2)(3t_1+2t_2)(\nr^x_1\nr^x_2+\nr^x_2\nr^x_1) \\
& \qquad +(3t_1+t_2)(2t_1+t_2)(3t_1+2t_2)(\nr^x_1\nr^x_2\nr^x_1+3\nr^x_2\nr^x_1\nr^x_2) \\
& \qquad +(3t_1+t_2)(2t_1+t_2)(3t_1+2t_2)(t_1+t_2)(\nr^x_1\nr^x_2\nr^x_1\nr^x_2+\nr^x_2\nr^x_1\nr^x_2\nr^x_1) \\
& \qquad +t_1(3t_1+t_2)(2t_1+t_2)(3t_1+2t_2)(t_1+t_2)\nr^x_1\nr^x_2\nr^x_1\nr^x_2\nr^x_1 \\
& \qquad +(3t_1+t_2)(2t_1+t_2)(3t_1+2t_2)(t_1+t_2)t_2\nr^x_2\nr^x_1\nr^x_2\nr^x_1\nr^x_2 \\
& \qquad +t_1(3t_1+t_2)(2t_1+t_2)(3t_1+2t_2)(t_1+t_2)t_2 \nr^x_1\nr^x_2\nr^x_1\nr^x_2\nr^x_1\nr^x_2 \\
& = h^x_2(t_2)h^x_1(t_1+t_2)h^x_2(3t_1+2t_2)h^x_1(2t_1+t_2)h^x_2(3t_1+t_2)h^x_1(t_1),
\end{align*}
which completes the proof.
\epf

Now we consider rank 2 Weyl groupoids with more than one point. Let $m=|\cX|$, and assume that the object change diagram has the form:
\begin{align}\label{eq:object change diagram}
& \vtxgpd^{x_1} \hspace{3pt} \text{\raisebox{3pt}{$\overset{2}{\rule{25pt}{0.5pt}}$}}
\hspace{3pt}\vtxgpd^{x_2} \hspace{3pt} \text{\raisebox{3pt}{$\overset{1}{\rule{25pt}{0.5pt}}$}}
\hspace{3pt} \vtxgpd^{x_3} \hspace{3pt} \text{\raisebox{3pt}{$\overset{2}{\rule{25pt}{0.5pt}}$}} \hspace{3pt} \dots
\hspace{3pt}\vtxgpd^{x_{m-1}} \hspace{3pt}
\text{\raisebox{3pt}{$\overset{j_{m-1}}{\rule{25pt}{0.5pt}}$}}
\hspace{3pt}\vtxgpd^{x_m} \hspace{3pt},
\end{align}
where $j_{2n}=1$, $j_{2n+1}=2$, $x_{n+1}=\rho_{j_n}(x_n)$, $n\in\Z$; notice that 
$x_k=x_{2rm+k}=x_{2rm-k+1}$ for $0\le k\le m$ and $r\in\Z$. 
Set also $i_{2n}=2$, $i_{2n+1}=1$.

Given an object $x$ of a Weyl groupoid of rank 2, we set
\begin{align*}
P^x_+&:= \prod_{n\alpha_1+m\alpha_2\in\Delta_+^x} (n \, t_1 + m \, t_2)\in\N_0[t_1,t_2].
\end{align*}
For each $\alpha\in\Delta_+^x$ we set also
\begin{align*}
P^x_+(\alpha)&:= \prod_{n\alpha_1+m\alpha_2\in\Delta_+^x\setminus\{\alpha\}} (n \, t_1 + m \, t_2)\in\N_0[t_1,t_2].
\end{align*}

\begin{lemma}\label{lemma:formula thm-several points}
For $x=x_1$,
\begin{align*}
h_{1}^{x_0}(t_1) & h_{2}^{x_1}(s_{1}^{x_1}\cdot t_2)\cdots h_{i_d}^{x_{d-1}}(s_{i_1}^{x_1}\cdots s_{i_{d-1}}\cdot t_{i_d})
\\ &= P_+^{x_1} \nr_{1}^{x_{0}}\nr_{2}^{x_{1}}\cdots \nr_{i_d}^{x_{d-1}}  + P_+^{x_1}(\alpha_1) \nr_{2}^{x_{1}}\cdots \nr_{i_d}^{x_{d-1}} \\
&=h_{2}^{x_1}(t_2)h_{1}^{x_2}(s_{2}^{x_1}\cdot t_1)\cdots h_{j_d}^{x_d}(s_{j_1}^{x_1}\cdots s_{j_{d-1}}\cdot t_{j_d}).
\end{align*}
For $x=x_k$, $1<k<m$ even,
\begin{align*}
h_{1}^{x_k}(t_1) & h_{2}^{x_{k+1}}(s_{1}^{x_k}\cdot t_2)\cdots h_{i_d}^{x_{k+d-1}}(s_{i_1}^{x_{k}}\cdots s_{i_{d-1}}\cdot t_{i_d})= P_+^{x_k} \nr_{1}^{x_{k}}\nr_{2}^{x_{k+1}}\cdots \nr_{i_d}^{x_{k+d-1}}\\
&=h_{2}^{x_k}(t_2)h_{1}^{x_{k-1}}(s_{2}^{x_k}\cdot t_1)\cdots h_{j_d}^{x_{k-d+1}}(s_{j_1}^{x_{k}}\cdots s_{j_{d-1}}\cdot t_{j_d}).
\end{align*}
\end{lemma}
\pf
Let $\beta_\ell=s_1^{x_k}s_2\dots s_{i_{\ell-1}}(\alpha_{i_\ell})$; in particular $\beta_1=\alpha_1$, $\beta_d=\alpha_2$, $\Delta_+^{a_k}=\{\beta_\ell:1\le \ell\le d\}$ and $\beta_1<\beta_2<\dots<\beta_d$ is the convex order corresponding to the reduced expression $w=s_1^{x_k}s_2\cdots s_{i_d}$. But $\beta_d<\beta_{d-1}<\dots<\beta_1$ is also a convex order so it corresponds to another reduced expression of $w$ by \cite{A-convex}: the unique possible one starting with $\delta_1=\alpha_2$ is $w=s_2^{x_k}s_1\dots s_{j_d}$. That is, if $\delta_\ell=s_2^{x_k}s_1\dots s_{j_{\ell-1}}(\alpha_{j_\ell})$, then $\delta_\ell=\beta_{d-\ell+1}$.

Notice that $m|d$ since $(\rho_1\rho_2)^d(x_k)=x_k$ for all $1\le k\le m$.

If $1<k<m$ and $k$ is even, then
\begin{align*}
h_{1}^{x_{k}}(t_1)&= t_1\nr_1^{x_{k}}, & h_{i_d}^{x_{k+d-1}}(s_{i_1}^{x_{k}}\cdots s_{i_{d-1}}\cdot t_{i_d})&=t_2\nr_{i_d}^{x_{k+d-1}}.
\end{align*}
Using the defining relations of $\cN$ in \eqref{eq:definign relations N} we see that
$$ h_{1}^{x_{k}}(t_1) h_{2}^{x_{k+1}}(s_{1}^{x_k}\cdot t_1)\cdots h_{i_d}^{x_{k+d-1}}(s_{i_1}^{x_{k}}\cdots s_{i_{d-1}}\cdot t_{i_d})
= P_+^{x_k} \nr_{1}^{x_{k}}\nr_{2}^{x_{k+1}}\cdots \nr_{i_d}^{x_{k+d-1}}, $$
since $\nr_i^{x_s}\nr_j^{x_t}=(\nr_i^{x_s}\er_i^{\rho_i(x_s)})(\er_j^{x_t}\nr_j^{x_t})=\delta_{\rho_i(x_s),x_t} \nr_i^{x_s}\nr_j^{x_t}$. Analogously,
$$ h_{2}^{x_{k}}(t_2)h_{1}^{x_{k-1}}(s_{2}^{x_{k}}\cdot t_2)\cdots h_{j_d}^{x_{k-d+1}}(s_{j_1}^{x_{k}}\cdots s_{j_{d-1}}\cdot t_{j_d})
= P_+^{x_k} \nr_{2}^{x_{k}}\nr_{1}^{x_{k-1}}\cdots \nr_{j_d}^{x_{k-d+1}}.$$
By \eqref{eq:indep reduced expression nil hecke}, $\nr_{1}^{x_{k}}\nr_{2}^{x_{k+1}}\cdots \nr_{i_d}^{x_{k+d-1}}=\nr_{2}^{x_{k}}\nr_{1}^{x_{k-1}}\cdots \nr_{j_d}^{x_{k-d+1}}$, so the equality between both expressions follows.

For the first equality, notice that $h_{i_d}^{x_{d-1}}(s_{i_1}^{x_{1}}\cdots s_{i_{d-1}}\cdot t_{i_d})=h_{i_d}^{x_{d-1}}(t_2)=t_2\nr_{i_d}^{x_{d-1}}$ since $m|d$, and $h_{1}^{x_1}(t_1)=\er^{x_1}+t_1\nr_1^{x_1}$, so the equality
$$ h_{1}^{x_0}(t_1) \cdots h_{i_d}^{x_{d-1}}(s_{i_1}^{x_1}\cdots s_{i_{d-1}}\cdot t_{i_d})=P_+^{x_1} \nr_{1}^{x_{0}}\nr_{2}^{x_{1}}\cdots \nr_{i_d}^{x_{d-1}}  + P_+^{x_1}(\alpha_1) \nr_{2}^{x_{1}}\cdots \nr_{i_d}^{x_{d-1}} $$
follows again by \eqref{eq:definign relations N}. Analogously
$$ h_{2}^{x_1}(t_2)\cdots h_{j_d}^{x_{d}}(s_{j_1}^{x_2}\cdots s_{j_{d-1}}\cdot t_{j_d})=P_+^{x_1} \nr_{2}^{x_{1}}\nr_{1}^{x_{2}}\cdots \nr_{j_d}^{x_{d}}  + P_+^{x_1}(\alpha_1) \nr_{2}^{x_{1}}\cdots \nr_{j_{d-1}}^{x_{d-1}}, $$
since $h_{j_d}^{x_{d}}(s_{j_1}^{x_1}\cdots s_{j_{d-1}}\cdot t_{j_d})=\er^{x_d}+t_1\nr_{j_d}^{x_d}$ and $h_{1}^{x_2}(t_2)=t_2\nr_2^{x_1}$. Finally the equality between both expressions follows again by \eqref{eq:indep reduced expression nil hecke}.
\epf

\begin{example}{\rm{
For $x_1$ as in \eqref{eq:cartan-graph standard A2}
\begin{align*}
h^{x_1}_1(t_1)h^{x_1}_2(t_1+t_2)h^{x_2}_1(t_2) & =(t_1+t_2)t_2\nr^{x_1}_2\nr^{x_2}_1+t_1(t_1+t_2)t_2\nr^{x_1}_1\nr^{x_1}_2\nr^{x_2}_1 \\
&= h^{x_1}_2(t_2)h^{x_2}_1(t_1+t_2)h^{x_3}_2(t_1),
\end{align*}
and for $x_2$,
\begin{align*}
h^{x_2}_1(t_1)h^{x_3}_2(t_1+t_2)h^{x_3}_1(t_2) & = t_1(t_1+t_2)t_2\nr^{x_2}_1\nr^{x_3}_2\nr^{x_3}_1 \\
&= h^{x_2}_2(t_2)h^{x_1}_1(t_1+t_2)h^{x_1}_2(t_1).
\end{align*}

Now for $x_1$ as in \eqref{eq:cartan-graph standard B2} we have that
\begin{align*}
h^{x_1}_1(t_1)&h^{x_1}_2(2t_1+t_2)h^{x_2}_1(t_1+t_2)h^{x_2}_2(t_2) = \\
& = (2t_1+t_2)(t_1+t_2)t_2\nr^{x_1}_2\nr^{x_2}_1\nr^{x_2}_2+t_1(2t_1+t_2)(t_1+t_2)t_2 \nr^{x_1}_1\nr^{x_1}_2\nr^{x_2}_1\nr^{x_2}_2 \\
& = h^{x_1}_2(t_2)h^{x_2}_1(t_1+t_2)h^{x_2}_2(2t_1+t_2)h^{x_1}_1(t_1).
\end{align*}
Notice the difference with these cases of standard Weyl groupoids with more than one point and those with a unique $x$ in Lemma \ref{lemma:formula thm-one point}; although they have the same sets of roots and the same finite Cartan matrices (of type $A_2$, $B_2$ respectively), the expressions of products in Lemmas \ref{lemma:formula thm-one point} and \ref{lemma:formula thm-several points} differ.
}}
\end{example}

\begin{example}{\rm{
We compute the expression in Lemma \ref{lemma:formula thm-several points} for the Weyl groupoid of the Nichols algebra \cite[Table 1, Row 10]{H-classif}. Its Cartan graph is
\begin{align*}
\setlength{\unitlength}{1mm}
\begin{picture}(40,3)(0,0)
\put(0,0){\circle*{1.5}}\put(1,2){{\small{$x_1$}}}
\put(4, 0){\line(1,0){11}}\put(9,1){{\small{$1$}}}
\put(17,0){\circle*{1.5}}\put(18,2){{\small{$x_2$}}}
\put(21, 0){\line(1,0){11}}\put(26,1){{\small{$2$}}}
\put(34,0){\circle*{1.5}}\put(35,2){{\small{$x_3$}}}
\end{picture}
\end{align*}
for the following Cartan matrices:
\begin{align*}
C_1&=\begin{bmatrix} 2 & \text{--}2 \\ \text{--}2 & 2 \end{bmatrix}, &  C_2&=\begin{bmatrix} 2 & \text{--}2 \\ \text{--}1 & 2 \end{bmatrix}, &
C_3&=\begin{bmatrix} 2 & \text{--}4 \\ \text{--}1 & 2 \end{bmatrix}.
\end{align*}
For $x_1$ we have that
\begin{align*}
h^{x_1}_1(t_1)&h^{x_2}_2(2t_1+t_2)h^{x_3}_1(t_1+t_2)h^{x_3}_2(2t_1+3t_2)h^{x_2}_1(t_1+2t_2)h^{x_1}_2(t_2) \\
& = P_+^{x_1} \,  \nr^{x_1}_1\nr^{x_2}_2\nr^{x_3}_1\nr^{x_3}_2\nr^{x_2}_1\nr^{x_1}_2 + P_+^{x_1}(\alpha_2) \,  \nr^{x_1}_1\nr^{x_2}_2\nr^{x_3}_1\nr^{x_3}_2\nr^{x_2}_1 \\
& = P_+^{x_1} \,  \nr^{x_1}_2\nr^{x_1}_1\nr^{x_2}_2\nr^{x_3}_1\nr^{x_3}_2\nr^{x_2}_1 + P_+^{x_1}(\alpha_2) \,  \nr^{x_1}_1\nr^{x_2}_2\nr^{x_3}_1\nr^{x_3}_2\nr^{x_2}_1 \\
& = h^{x_1}_2(t_2)h^{x_1}_1(t_1+2t_2)h^{x_2}_2(2t_1+3t_2)h^{x_3}_1(t_1+t_2)h^{x_3}_2(2t_1+t_2)h^{x_2}_1(t_1).
\end{align*}
For $x_2$ we compute
\begin{align*}
h^{x_2}_1&(t_1)h^{x_1}_2(2t_1+t_2)h^{x_1}_1(3t_1+2t_2)h^{x_2}_2(4t_1+3t_2)h^{x_3}_1(t_1+t_2)h^{x_3}_2(t_2) \\
& = P_+^{x_2} \,  \nr^{x_2}_1\nr^{x_1}_2\nr^{x_1}_1\nr^{x_2}_2\nr^{x_3}_1\nr^{x_3}_2 = P_+^{x_2} \,  \nr^{x_2}_2\nr^{x_3}_1\nr^{x_3}_2\nr^{x_2}_1\nr^{x_1}_2\nr^{x_1}_1 \\
& = h^{x_2}_2(t_2)h^{x_3}_1(t_1+t_2)h^{x_3}_2(4t_1+3t_2)h^{x_2}_1(3t_1+2t_2)h^{x_1}_2(2t_1+t_2)h^{x_1}_1(t_1).
\end{align*}
Finally for $x_3$ we have that
\begin{align*}
h^{x_3}_1(t_1)&h^{x_3}_2(4t_1+t_2)h^{x_2}_1(3t_1+t_2)h^{x_1}_2(2t_1+t_2)h^{x_1}_1(t_1+t_2)h^{x_2}_2(t_2) \\
& = P_+^{x_3} \,  \nr^{x_3}_1\nr^{x_3}_2\nr^{x_2}_1\nr^{x_1}_2\nr^{x_1}_1\nr^{x_2}_2 + P_+^{x_3}(\alpha_1) \,  \nr^{x_3}_2\nr^{x_2}_1\nr^{x_1}_2\nr^{x_1}_1\nr^{x_2}_2 \\
& = P_+^{x_3} \,  \nr^{x_3}_2\nr^{x_2}_1\nr^{x_1}_2\nr^{x_1}_1\nr^{x_2}_2\nr^{x_3}_1 + P_+^{x_3}(\alpha_1) \,  \nr^{x_3}_2\nr^{x_2}_1\nr^{x_1}_2\nr^{x_1}_1\nr^{x_2}_2 \\
& = h^{x_3}_2(t_2)h^{x_2}_1(t_1+t_2)h^{x_1}_2(2t_1+t_2)h^{x_1}_1(3t_1+t_2)h^{x_2}_2(4t_1+t_2)h^{x_3}_1(t_1).
\end{align*}
}}
\end{example}

\begin{prop}\label{prop:toolone}
Fix $i\ne j$, $x\in \cX$. For $n\in\N_0$, set $m=m_{ij}^x$,
\begin{align}\label{eq:i,j alternate}
i_{2n}&=j, &  i_{2n+1}&=i, & x_n&=\rho_{i_{n-1}}\cdots \rho_{i_1}(x),\\
j_{2n}&=i, &  j_{2n+1}&=j, & y_n&=\rho_{i_{n-1}}\cdots \rho_{j_1}(x).\notag
\end{align}
Then
\begin{align*}
h_{i_1}^{x_1}(t_{i_1}) & h_{i_2}^{x_2}(s_{i_1}^{x_1}\cdot t_{i_2})\cdots h_{i_m}^{x_m}(s_{i_1}^{x_1}\cdots s_{i_{m-1}}\cdot t_{i_m})= \\
&=h_{j_1}^{y_1}(t_{j_1})h_{j_2}^{y_2}(s_{j_1}^{y_1}\cdot t_{j_2})\cdots h_{j_m}^{y_m}(s_{j_1}^{y_1}\cdots s_{j_{m-1}}\cdot t_{j_m}).
\end{align*}
\end{prop}
\pf
If $\rho_i(x)=\rho_j(x)=x$, then it follows by Lemma \ref{lemma:formula thm-one point}. If either $\rho_i(x)\neq x$ or $\rho_j(x)\neq x$, and the object change is not a cycle, then it has the shape \eqref{eq:object change diagram} and the equality follows by Lemma \ref{lemma:formula thm-several points}.
Otherwise the object change is a cycle with $d$ objects, where $d|m$. Then $h_{i_\ell}^{x_\ell}(t)=t\nr_{i_\ell}^{x_\ell}$ for all $\ell$ so
\begin{align*}
h_{i_1}^{x_1}(t_{i_1}) & h_{i_2}^{x_2}(s_{i_1}^{x_1}\cdot t_{i_2})\cdots h_{i_m}^{x_m}(s_{i_1}^{x_1}\cdots s_{i_{m-1}}\cdot t_{i_m})
= P_+^{x_1} \nr_{1}^{x_1}\nr_{2}^{x_2}\cdots \nr_{i_m}^{x_m}\\
&=h_{j_1}^{y_1}(t_{j_1})h_{j_2}^{y_2}(s_{j_1}^{y_1}\cdot t_{j_2})\cdots h_{j_m}^{y_m}(s_{j_1}^{y_1}\cdots s_{j_{m-1}}\cdot t_{j_m}),
\end{align*}
where last equality follows by \eqref{eq:indep reduced expression nil hecke}.
\epf

Hence we obtain a representation of the groupoid $\widetilde{\cW}$, a key step towards our main result, Theorem \ref{thm:subsequences bruhat order}.

\begin{theorem}
There exists a well-defined map $\Psi:\widetilde{\cW}\to\End(\bK \cN)$ such that
\begin{align}\label{eq:def Psi}
\Psi(\st_{i_1}^x \cdots \st_{i_m})(f \, T_w)&:=
h^{x_1}_{i_1}(\beta_1)h^{x_2}_{i_2}(\beta_2)\cdots h^{x_m}_{i_1}(\beta_m) ( s_{i_1}^{x_1} \cdots s_{i_m}^{x_m} \cdot f) \, T_w,
\end{align}
for each $f\in\bK$, $w\in\cW$, where $x_j=\rho_{i_{j-1}} \cdots \rho_{i_1}(x)$, $\beta_j=s_{i_1}^x \cdots s_{i_{j-1}}(\al_{i_j})$.
\end{theorem}
\pf
We claim that there exists $\Psi:\widetilde{\cW}\to\End(\bK \cN)$ such that
\begin{align*}
\Psi(\id_x \st_i)(f \, T_w)&:= h^a_{i}(\alpha_i)( s_{i}^a \cdot f) \, T_w, & f\in\bK,&\,w\in\cW,
\end{align*}
from which the statement follows. We have to prove that it is well defined. Fix $i\ne j$, $x\in \cX$. Set $m=m_{ij}^x$, and the same elements defined in \eqref{eq:i,j alternate}. It is enough to prove that
$$ \Psi(\st_{i_1}^x \cdots \st_{i_d})= \Psi(\st_{j_1}^x \cdots \st_{j_d}).$$
But this equality follows by Proposition \ref{prop:toolone} and \cite[Lemma 5]{HY-weylgpd}, since
\begin{align*}
\Psi(\st_{i_1}^x \cdots \st_{i_d})(f \, T_w)&:= h_{i_1}^{x_1}(t_{i_1})h_{i_2}^{x_2}(s_{i_1}^{x_1}\cdot t_{i_2})\cdots h_{i_d}^{x_d}(s_{i_1}^{x_1}\cdots s_{i_{d-1}}\cdot t_{i_d})
\\ & \qquad ( s_{i_1}^{x_1} \cdots s_{i_d} \cdot f) \, T_w,
\end{align*}
and we have the corresponding formula for $\Psi(\st_{j_1}^x \cdots \st_{j_d})(f \, T_w)$.
\epf

\section{Bruhat order for Weyl groupoids}\label{section:Bruhat order}

Now we are ready to define the Bruhat order for Weyl groupoids.

\begin{definition}
Let $(\cX, \rho, \cC)$ be a semi-Cartan graph and $\cR$ a GRS of $(\cX, \rho, \cC)$.

Let $x,y\in\cX$. For $w \in \cW(y,x)$
with $s_{i_1}^{x_1} s_{i_2}^{x_2}\cdots s_{i_r}^{x_r}$ being a reduced expression of $w$. 
We say that a subsequence $a(1),\ldots,a(t)$ of $i_1\ldots,i_r$
is {\it{$(y,x)$-good}} if $x_{z+1}=x_z$ for $z\in\{1,\ldots,r\}\setminus\{a(1),\ldots,a(t)\}$.
Given $u,w \in \cW(y,x)$, we say that $u<w$ if there exists a reduced expression
$w=s_{i_1}^{x_1} s_{i_2}^{x_2}\cdots s_{i_r}^{x_r}$ and a subsequence $a(1),\ldots,a(t)$ 
of $i_1\ldots,i_r$ such that $u=s_{a(1)}^x s_{a(2)}\cdots s_{a(t)}$ is a reduced expression of $u$. It defines a partial order on $\cW(y,x)$ called the \emph{Bruhat order}.
\end{definition}


Now we are ready to prove that Bruhat order is well-defined.

\begin{theorem}\label{thm:subsequences bruhat order}
Let $w\in\cW(y,x)$, $x,y\in\cX$. Assume that $r:=\ell(w)>0$ and
$$ w=\sigma^{x}_{i_1}\cdots \sigma_{i_r} = \sigma^{x}_{j_1}\cdots s_{j_r}.$$
are two reduced expressions of $w$.

Then for any {\it{$(y,x)$-good}} sub-sequence $a(1),\ldots,a(t)$ of $i_1,\ldots,i_r$ such that
\begin{align*}
u&:=\sigma^{x}_{{a(1)}}\cdots \sigma_{{a(t)}}, & \ell(u)&=t, \, u\in\cW(y,x),
\end{align*}
there exists a {\it{$(y,x)$-good}} sub-sequence $b(1),\ldots,b(t)$ of $j_1,\ldots,j_r$ such that
$$ u=\sigma^{x}_{{b(1)}}\cdots \sigma_{{b(t)}}. $$
\end{theorem}
\pf
Fix $w\in\cW(y,x)$ such that $r:=\ell(w)>0$, and a reduced expression $ w=\sigma^{x}_{i_1}\cdots \sigma_{i_r}$.
Let $\Theta(x;i_1,\ldots,i_r)$ be the set of elements
\begin{align*}
u&:=\sigma^{x}_{{a(1)}}\cdots \sigma_{{a(t)}}, & \mbox{satisfies }\ell(u)&=t, \, u\in\cW(y,x),
\end{align*}
for some $(y,x)$-good sub-sequence $a(1),\ldots,a(t)$ of $i_1,\ldots,i_r$.

Let $u\in\cW$, $x\in\cX$, $i\in I$, $m=\ell(u)$. From Theorem \ref{thm:basis nil Hecke} we see that $n_i^x T_u=0$ if $\ell(\sigma_i^x u)<m+1$, and $n_i^x T_u=T_{\sigma_i^x u}$ if $\ell(\sigma_i^x u)=m+1$. Note that
\begin{align*}
\Psi(\st_{i_1}^x \cdots \st_{i_m})(1)&=h^{x_1}_{i_1}(\beta_1)h^{x_2}_{i_2}(\beta_2)\cdots h^{x_m}_{i_1}(\beta_m),
\end{align*}
and all the $\beta_j$ are positive roots. Then there exist non-zero polynomials $f_u\in\Z[t_i|i\in I]$, $u\in\Theta(X_1;i_1,\ldots,i_r)$, such that
$$ \Psi(\st_{i_1}^x \cdots \st_{i_m})(1) = \sum_{u\in \Theta(x;i_1,\ldots,i_r)} f_u \, T_u. $$
Moreover all the polynomials $f_u$ have non-negative coefficients.

Therefore the set $\Theta(X_1;i_1,\ldots,i_r)$ is identified with the subset of $\cW$ such that $T_u$ appears with non-zero coefficient in
$\Psi( \st_{i_1}^x \cdots \st_{i_m})(1)$. But $\Psi( \st_{i_1}^x \cdots \st_{i_m})(1)$ does not depend on the reduced expression of $w$, so the statement follows.
\epf

\end{document}